\documentclass[11pt]{article}
\usepackage{amsfonts,latexsym,amsmath,amscd,geometry}
\usepackage{amssymb}
\usepackage{latexsym}

\newcommand \nc{\newcommand}
\newtheorem{theorem}{Theorem}[section]
\newtheorem{lemma}[theorem]{Lemma}

\newtheorem{definition}[theorem]{Definition}

\newtheorem{remark}[theorem]{Remark}

\nc{\ba}{\begin{array}}\nc{\ea}{\end{array}}
\nc{\be}{\begin{eqnarray}}\nc{\ee}{\end{eqnarray}}
\nc{\beq}{\begin{equation}}\nc{\eeq}{\end{equation}}
\nc{\bex}{\begin{eqnarray*}}\nc{\eex}{\end{eqnarray*}}
\nc{\btm}{\begin{theorem}} \nc{\etm}{\end{theorem}}
\nc{\blm}{\begin{lemma}} \nc{\elm}{\end{lemma}}
\nc{\R}{\mathbb{R}} \nc{\va}{\varepsilon} \nc{\ls}{\limits}

\def\pf{\noindent{\bf Proof.\quad}}
\def\les{\lesssim}

\newcommand \qed {\hfill $\Box$}

\begin{document}
\title{Boundary regularity of stationary biharmonic maps}
\author{Huajun Gong\footnote{School of Mathematics, Fudan University, Shanghai 200433, P.R. China}
\quad Tobias Lamm\footnote{Institut f\"ur  Mathematik, Goethe-Universit\"at Frankfurt, Robert-Mayer-Str.1060054, Frankfurt, Germany}\quad
Changyou Wang\footnote{Department of Mathematics, University of Kentucky, Lexington, KY 40506, USA}}
\date{}
\maketitle

\begin{abstract} We consider the Dirichlet problem for stationary biharmonic maps $u$ from a bounded, smooth domain
$\Omega\subset\mathbb R^n$ ($n\ge 5$) to a compact, smooth Riemannian manifold $N\subset\mathbb R^l$
without boundary. For any smooth boundary data, we show that if, in addition, $u$ satisfies a certain boundary monotonicity inequality, then there exists a closed subset $\Sigma\subset\overline{\Omega}$,
with $H^{n-4}(\Sigma)=0$, such that $u\in C^\infty(\overline\Omega\setminus\Sigma, N)$.
\end{abstract}

\setcounter{section}{0} \setcounter{equation}{0}
\section{Introduction}

This is a continuation of our previous study in \cite{LW}. Here we consider the Dirichlet problem
for (extrinsic) biharmonic maps into Riemannian manifolds in dimension at least $5$ and address
the issue of boundary regularity for a class of stationary biharmonic maps.

For $n\ge 5$, let $\Omega\subset\mathbb R^n$ be a bounded, smooth domain and
$(N, h)\subset\mathbb R^L$ be a $l$-dimensional, compact $C^3$-Riemannian manifold
 with $\partial N=\emptyset$.   For $k\ge 1, 1\le p<+\infty$, define
the Sobolev space
$$W^{k,p}(\Omega, N)=\Big\{v\in W^{k,p}(\Omega, \mathbb R^L): \ v(x)\in N
\ {\rm{for\ a.e.\ }} x\in\Omega\Big\}.$$
Recall that an extrinsic biharmonic map $u\in W^{2,2}(\Omega, N)$ is defined to be a critical
point of the Hessian energy functional:
$$E_2(v)=\int_\Omega |\Delta v|^2, \ v\in W^{2,2}(\Omega, N).$$

If we denote by $\mathbb P(y):\mathbb R^L\to T_y N$, $y\in N$, the orthogonal projection map, then
the second fundamental form of $B$ is defined by 
$$\mathbb B(y)(X, Y)=-D_X\mathbb P(y)(Y), \ \forall X, Y\in T_y N.$$
It is standard (cf. \cite{W1} Proposition 2.1) that an extrinsic biharmonic map $u\in W^{2,2}(\Omega, N)$
is a weak solution to the biharmonic map equation:
\begin{equation}\label{biharm1}
\Delta^2 u=\Delta(\mathbb B(u)(\nabla u,\nabla u))
+2\nabla\cdot\langle\Delta u, \nabla(\mathbb P(u))\rangle
-\langle\Delta (\mathbb P(u)),\Delta u\rangle,\end{equation}
or equivalently
\begin{equation}
\Delta^2 u\perp T_u N. \label{biharm2}
\end{equation}

Notice that the biharmonic map equation (\ref{biharm1}) is a 4th order elliptic system with super-critical
nonlinearity. It is a very natural and interesting question to study its regularity. The study was first initiated by Chang, Wang, and Yang \cite{CWY}. In particular, they proved that when $N=S^{L-1}\subset\mathbb R^L$ is the unit sphere, then any $W^{2,2}$-biharmonic map is smooth in dimension $4$, and smooth away from
a closed set of $(n-4)$-dimensional Hausdorff measure zero for $n\ge 5$ provided that it is, in addition, stationary. The main theorem in \cite {CWY} was subsequently extended to any smooth Riemannian manifold $N$ by the third author in \cite{W,W1, W2}, and different proofs were later given by Lamm-Rivier\'e \cite{LR}
for $n=4$ and Struwe \cite{SM} for $n\ge 5$, see also Strzelecki \cite{Str} for some generalizations. 

There have been many important works on the regularity of stationary harmonic maps, originally due to 
H\'elein \cite{H}, Evans \cite{E}, and Bethuel \cite{B} (see Rivier\'e \cite{RT}, and Rivier\'e-Struwe \cite{RS}
for important new approaches and improvements). A crucial property for stationary harmonic maps is the well-known energy monotonicity formula (see Price \cite{P}). The notion of stationary biharmonic maps was motived by the notion of stationary harmonic maps. More precisely, $u\in W^{2,2}(\Omega, N)$ is called
a stationary biharmonic map if it is, in addition, a critical point with respect to the domain variations:
\begin{equation}\label{stationary}
\frac{d}{dt}\Big|_{t=0}\int_{\Omega} |\Delta u_t|^2=0, \ u_t(x)=u(x+t Y(x)), 
\ Y\in C_0^\infty(\Omega,\mathbb R^n).
\end{equation}
It has been derived by \cite{CWY} and Angelsberg \cite{A} that stationary biharmonic maps enjoy the 
following interior monotonicity inequality: for $x\in\Omega$ and $0<r\le R<{\rm{dist}}(x,\partial\Omega)$,
\begin{equation}\label{interior_mono}
R^{4-n}\int_{B_R(x)}|\Delta u|^2-r^{4-n}\int_{B_r(x)}|\Delta u|^2=A_1+A_2,
\end{equation}
where
\begin{eqnarray*} A_1&=&4\int_{B_R(x)\setminus B_r(x)}
\left(\frac{|u_i+(y-x)^ju_{ij}|^2}{|y-x|^{n-2}}+(n-2)\frac{|(y-x)^i u_i|^2}{|y-x|^n}\right),\\
A_2&=&2\int_{\partial (B_R(x)\setminus B_r(x))}
\left(-\frac{(y-x)^i u_j u_{ij}}{|y-x|^{n-3}}+2\frac{|(y-x)^i u_i|^2}{|y-x|^{n-1}}
-2\frac{|\nabla u|^2}{|y-x|^{n-3}}\right).
\end{eqnarray*}
Such an interior monotonicity formula plays a critical role in the partial interior regularity theory, that was mentioned above, for stationary biharmonic maps.

It is also a natural and interesting question to address possible boundary regularity for biharmonic maps,
associated with smooth Dirichlet boundary data. More precisely, let 
$\phi\in C^\infty(\Omega_\delta, N)$ be given for some $\delta>0$, where 
$\Omega_\delta=\Big\{x\in\overline\Omega: {\rm{dist}}(x,\partial\Omega)<\delta\Big\}$. Consider that
$u\in W^{2,2}(\Omega, N)$ is a biharmonic map satisfying
\begin{equation}\label{bdry_cond0}
\Big(u,\ \frac{\partial u}{\partial\nu}\Big)\Big|_{\partial\Omega}
=\Big(\phi, \ \frac{\partial\phi}{\partial\nu}\Big),
\end{equation}
where $\nu$ is the unit outward normal of $\partial\Omega$.

For dimension of $\Omega$, $n=4$,  the complete boundary smoothness of biharmonic maps
has been proved by Ku \cite{K} for $N=S^{L-1}$ and Lamm-Wang \cite{LW} for any compact
Riemannian manifold $N$\footnote{In fact, the complete interior and boundary smoothness of 
$k$-polyharmonic maps in $W^{k,2}(\mathbb R^{2k}, N)$ has been proved
by Gastel-Scheven \cite{GS} and \cite{LW} respectively.}.  For dimensions $n\ge 5$,  as in the
interior case, it seems necessary to require a boundary monotonicity inequality analogous to
(\ref{interior_mono}) in order to obtain possible boundary regularity.  Here we introduce
a notion of {\it globally stationary biharmonic map} which enjoys a boundary monotonicity inequality.

\begin{definition} A stationary biharmonic map $u\in W^{2,2}(\Omega, N)$ associated with (\ref{bdry_cond0})
is called a globally stationary biharmonic map, if there exist $R_0>0$ and $C$, depending only on
$n, \partial\Omega, \delta, \|\phi\|_{C^4(\Omega_\delta)}$, such that for any $x_0\in\partial\Omega$
and $0<r\le R\le R_0$, there holds
\begin{eqnarray}\label{bdry_mono}
&&r^{4-n}\int_{\Omega\cap B_r(x_0)}|\nabla^2(u-\phi)|^2+A(r,R)+e^{Cr} B(r)+e^{cr}C(r)\nonumber\\
&\leq& CR e^{CR}+e^{CR} R^{4-n}\int_{\Omega\cap B_R(x_0)}|\nabla^2(u-\phi)|^2+e^{CR} B(R)
+e^{CR} C(R),
\end{eqnarray}
where
$$A(r,R):=4\int_{\Omega\cap(B_R(x_0)\setminus B_r(x_0))}
\left(\frac{|(u-\phi)_j+(y-x_0)^i (u-\phi)_{ij}|^2}{|y-x_0|^{n-2}}+(n-2)\frac{|(y-x_0)^i(u-\phi)_i|^2}{|y-x_0|^n}\right),$$
and for $\rho=r$ and $R$,
$$
B(\rho):=2\int_{\partial B_\rho(x_0)\cap\Omega}\left(2\frac{|(y-x_0)^k(u-\phi)_k|^2}{\rho^{n-1}}-\frac{(u-\phi)_i (u-\phi)_{ij} (y-x_0)^j}{\rho^{n-3}}
-2\frac{|\nabla (u-\phi)|^2}{\rho^{n-3}}\right),
$$
$$
C(\rho):=2\rho^{4-n}\int_{\partial B_\rho(x_0)\cap\Omega}
\Big[\langle \Delta(u-\phi),\frac{\partial}{\partial r} (u-\phi)\rangle
-\langle\nabla(u-\phi), \frac{\partial}{\partial r}(\nabla(u-\phi))\rangle\Big].
$$
\end{definition}

Now we state our main theorem.

\begin{theorem}\label{bdry_reg} For $n\ge 5$,\ $\phi\in C^\infty(\Omega_\delta,N)$ for some
$\delta>0$, if $u\in W^{2,2}(\Omega, N)$ associated with (\ref{bdry_cond0}) is a globally stationary biharmonic map, then there exists a closed subset $\mathcal S(u)\subset\overline\Omega$,
with $H^{n-4}(\mathcal S(u))=0$, such that $u\in C^\infty\Big(\overline\Omega\setminus\mathcal S(u), N\Big)$.
\end{theorem}

\begin{remark}  Notice that the interior monotonicity inequality (\ref{interior_mono}) is a consequence
of the stationarity identity (\ref{stationary}).  However, it is unclear whether the boundary monotonicity
inequality (\ref{bdry_mono}) can be deduced from (\ref{stationary}).  On the other hand, in \S 2
below, we will employ a Pohazev type argument to show that (\ref{bdry_mono}) holds for any 
sufficiently regular biharmonic map, e.g. $u\in W^{4,2}(\Omega_\delta, N)$. Thus it seems to be a necessary condition for boundary regularity.
\end{remark}

In contrast with the standard reflection argument to obtain boundary regularity of harmonic maps (see, e.g. Wang \cite{W3}), it seems impossible to obtain boundary regularity of biharmonic maps or general 4th order
elliptic systems via boundary reflection methods. To overcome this type of difficulty,
we will employ an estimate for the Green function of $\Delta^2$ on the modified upper half ball
to prove a boundary decay lemma under a smallness condition in suitable Morrey spaces.
The overall scheme includes suitable adaptions and extensions of: 
(i) the rewriting of biharmonic map equation (\ref{biharm1}) by
\cite{LR} and \cite{SM}, and the Coloumb gauge construction by \cite{RS} and \cite{SM}; 
and (ii) estimate of Riesz potentials among Morrey spaces by \cite{W1,W2}. To obtain
the Morrey bound (\ref{bdry_small}), we modify the argument by \cite{SM} appendix B.

The paper is organized as follows. In \S2, we derive the boundary monotonicity inequality for
$W^{4,2}$-biharmonic maps, which may have its own interests. In \S3, we establish a boundary
decay lemma under the smallness condition (\ref{bdry_small}). In \S4, we establish (\ref{bdry_small})
at $H^{n-4}$ a.e. $x\in\partial\Omega$.

Throughout this paper, for two quantities $A$ and $B$,
we will denote by $A\lesssim B$ if there exists a constant $C>0$,
depending only on $n, \Omega$, and $\phi$, such that $A\le CB$.

\setcounter{section}{1} \setcounter{equation}{0}
\section{Derivation of boundary monotonicity inequality for  regular biharmonic maps}
This section is devoted to the derivation of  the boundary Hessian energy monotonicity
inequality (\ref{bdry_mono})  for $W^{4,2}(\Omega_\delta)$-extrinsic biharmonic maps.

For $n\ge 5$, set  $\mathbb R^n_+=\Big\{x=(x',x^n)\in\mathbb R^n: x^n\ge 0\Big\}$.
For $R>0$, set
$$B_R=\Big\{x\in\mathbb R^n:\ |x|\le R\Big\} 
\ \ {\rm{and}}\ \ B_R^+=B_R\cap \mathbb R^n_+.$$ 
Denote by 
$$T_R=\partial B_R^+\cap\Big\{x=(x',x^n):\ x^n=0\Big\}\ \ {\rm{and}}\ \ \ S_R^+=\partial B_R\cap\Big\{x=(x',x^n)
:\ x^n>0\Big\}$$
the flat part and curved part of $\partial B_R^+$ respectively.

For simplicity, we will derive (\ref{bdry_mono}) for the case that $x_0=0$, $R_0=2$, and $\Omega\cap B_{R_0}(x_0)=B_2^+$.
Hence we may assume that $\phi\in C^4(\overline{B_2^+}, N)$ and $u\in W^{4,2}(B_2^+,N)$ is a  
biharmonic map satisfying
\begin{equation} \label{bdry_cond}
u\Big|_{T_2}=\phi\Big|_{T_2},\  \ {\partial u\over\partial x_n}\Big|_{T_2}
={\partial\phi\over\partial x_n}\Big|_{T_2}.
\end{equation}

We want to show that $u$ satisfies (\ref{bdry_mono}) for $0<r\le R\le 1$.
For $0<r\le 1$,  set $v_r(x)=v(rx)$ for $x\in B_2^+$.
Direct compuation using integration by parts implies 
\begin{eqnarray}
&&\frac{d}{dr}\left(r^{4-n}\int_{B_r^+}|\Delta (u-\phi)|^2\right)
=\frac{d}{dr}\int_{B_1^+}\Big|\Delta (u-\phi)_r\Big|^2\nonumber\\
&=&2\int_{B_1^+}\langle\Delta (u-\phi)_r(x),\ \Delta\left(x\cdot\nabla (u-\phi)(rx)\right)\rangle\nonumber\\
&=&2r^{3-n}\int_{B_r^+}\left\langle\Delta (u-\phi),\ \Delta\left(x\cdot\nabla (u-\phi)\right)\right\rangle
\nonumber\\
&=&2r^{3-n}\int_{B_r^+}\nabla\cdot\langle \Delta (u-\phi),\ \nabla\left(x\cdot\nabla(u-\phi)\right)\rangle
-2r^{3-n}\int_{B_r^+}\left\langle \nabla\Delta (u-\phi),\ \nabla\left(x\cdot\nabla(u-\phi)\right)\right\rangle\nonumber\\
&=& 2r^{3-n}\int_{S_r^+}\left\langle\Delta (u-\phi),\ {\partial\over\partial r}\left(r{\partial (u-\phi)\over\partial r}\right)\right\rangle\nonumber\\
&-&2r^{3-n}\int_{B_r^+}\left\langle \nabla\Delta (u-\phi),\ \nabla\left(x\cdot\nabla(u-\phi)\right)\right\rangle,\label{radial1}
\end{eqnarray}
where we have used  the fact 
\begin{eqnarray*}&&\int_{T_r}\left\langle\Delta (u-\phi), \ {\partial\over\partial x_n}\left(x\cdot\nabla(u-\phi)\right)\right\rangle\\
&=&\int_{T_r}\left\langle\Delta (u-\phi),\ \frac{\partial}{\partial x_n}\left(u-\phi\right)+\sum_{j=1}^n x^j 
\frac{\partial^2}{\partial x_j\partial x_n}\left(u-\phi\right)\right\rangle =0,
\end{eqnarray*}
since 
$$x^n=\frac{\partial}{\partial x_n}\left(u-\phi\right)=\frac{\partial^2}{\partial x_j\partial x_n}\left(u-\phi\right)=0
\ \ {\rm{for}}\ \ 1\le j\le n-1\ \  \ {\rm{on}}\ \  T_r.$$

Now we calculate, using integration by parts,
\begin{eqnarray}
&&\int_{B_r^+}\left\langle \nabla\Delta (u-\phi),\ \nabla\left(x\cdot\nabla(u-\phi)\right)\right\rangle\nonumber\\
&=&\int_{B_r^+}\nabla\cdot\left\langle\nabla\Delta (u-\phi), \ x\cdot\nabla\left(u-\phi\right)\right\rangle
-\int_{B_r^+}\left\langle \Delta^2 (u-\phi), \ x\cdot\nabla\left(u-\phi\right)\right\rangle\nonumber\\
&=&\int_{S_r^+}\left\langle{\partial\over\partial r}(\Delta (u-\phi)),\  x\cdot\nabla\left(u-\phi\right)\right\rangle
\nonumber\\
&+&\int_{B_r^+}\left\langle\Delta^2 u,\ x\cdot\nabla\phi\right\rangle
+\int_{B_r^+}\langle\Delta^2\phi, \ x\cdot\nabla(u-\phi)\rangle,\label{radial2}
\end{eqnarray}
where we have used  in the last step that 
$$ \Big\langle\Delta^2 u, \ x\cdot\nabla u\Big\rangle =0 \ \ {\rm{on}}\ \  B_r^+ \ 
\Big(\Leftarrow \Delta^2 u\perp T_u N
\ {\rm{and}}\ x\cdot\nabla u\in T_u N\Big), $$
and
$$\left\langle{\partial\over\partial x_n}(\Delta u), \ x\cdot\nabla\left(u-\phi\right)\right\rangle=0\ \ \hbox{ on }\ \ T_r.$$
Substituting (\ref{radial2}) into (\ref{radial1}), we obtain
\begin{eqnarray}
&&{d\over dr}\left(r^{4-n}\int_{B_r^+}|\Delta (u-\phi)|^2\right)\nonumber\\
&=&2r^{3-n}\int_{S_r^+}\Big[\left\langle\Delta (u-\phi),\ {\partial\over \partial r}\left(x\cdot\nabla (u-\phi)\right)\right\rangle-\left\langle{\partial\over \partial r}(\Delta (u-\phi)),\ x\cdot\nabla(u-\phi)\right\rangle\Big]\nonumber\\
&-&2r^{3-n}\int_{B_r^+}\left\langle\Delta^2\phi,\ x\cdot\nabla (u-\phi)\right\rangle
-2r^{3-n}\int_{B_r^+}\left\langle \Delta^2 u,\ x\cdot\nabla\phi\right\rangle. \label{radial3}
\end{eqnarray}

It is easy to see that by Poincar\'e's inequality,
the third term in the right hand side of (\ref{radial3}) can be estimated by
\begin{eqnarray} \label{thirdterm-est}
\Big|2r^{3-n}\int_{B_r^+}\left\langle\Delta^2 \phi,\ x\cdot\nabla(u- \phi)\right\rangle\Big|
&\lesssim& r^{2-n}\int_{B_r^+}|\nabla (u-\phi)|^2 +r^4\nonumber\\
&\lesssim& r^{4-n}\int_{B_r^+}|\nabla^2(u-\phi)|^2+r^4.
\end{eqnarray}
Now we estimate the last term in the right hand side of (\ref{radial3}).  Applying (\ref{biharm1}), we have
\begin{eqnarray*}
\int_{B_r^+}\Big\langle\Delta^2 u, \ x\cdot\nabla\phi\Big\rangle
&=&\int_{B_r^+}\Big\langle\Delta\Big(\mathbb B(u)(\nabla u,\nabla u)\Big),\ x\cdot\nabla\phi\Big\rangle\\
&+&2\int_{B_r^+}\Big\langle\nabla\cdot\Big\langle\Delta u,\ \nabla(\mathbb P(u))\Big\rangle,\ x\cdot\nabla\phi\Big\rangle\\
&-&\int_{B_r^+}\Big\langle\Big\langle\Delta(\mathbb P(u)),\ \Delta u\Big\rangle, \ x\cdot\nabla\phi\Big\rangle\\
&=&I+II+III.\end{eqnarray*}
Now we estimate $I, II, III$ as follows.
\begin{eqnarray*}
I&=&\int_{B_r^+}\nabla\cdot\Big\langle\nabla\Big(\mathbb B(u)(\nabla u,\nabla u)\Big),\ x\cdot\nabla\phi\Big\rangle
-\int_{B_r^+}\Big\langle\nabla\Big(\mathbb B(u)(\nabla u,\nabla u)\Big),\ \nabla(x\cdot\nabla\phi)\Big\rangle\\
&=&I_1+I_2.\end{eqnarray*}
It is easy to see that
\begin{eqnarray}\label{I2-est}
\Big|I_2\Big|&\lesssim &\int_{B_r^+}\Big(|\nabla u|^3+|\nabla u||\nabla^2 u|\Big)
\Big(|\nabla \phi|+|x||\nabla^2\phi|\Big)\nonumber\\
&\lesssim & \int_{B_r^+}\Big(|\nabla u|^3+|\nabla u||\nabla^2 u|\Big),
\end{eqnarray}
while
\begin{eqnarray*}
I_1&=&\int_{S_r^+}\Big\langle{\partial\over\partial r}\Big(\mathbb B(u)(\nabla u,\nabla u)\Big),
\ x\cdot\nabla\phi\Big\rangle-\int_{T_r}\Big\langle {\partial\over\partial x_n}\Big(\mathbb B(u)(\nabla u,\nabla u)\Big),
\ x\cdot\nabla \phi\Big\rangle\\
&=&I_{1a}+I_{1b}\end{eqnarray*}
To estimate $I_{1b}$, observe that since $\Big\langle \mathbb B(u)(\nabla u,\nabla u),\ x\cdot\nabla u\Big\rangle=0$,
we have
\begin{eqnarray*}
I_{1b}&=&-\int_{T_r}\Big\langle {\partial\over\partial x_n}\Big(\mathbb B(u)(\nabla u,\nabla u)\Big),
\ x\cdot\nabla u\Big\rangle\\
&=&\int_{T_r}\Big\langle \mathbb B(u)(\nabla u,\nabla u),
\ \frac{\partial}{\partial x_n}(x\cdot\nabla u)\Big\rangle.
\end{eqnarray*}
By using the boundary condition (\ref{bdry_cond}), we have
$$\frac{\partial}{\partial x_n}(x\cdot\nabla u)=\sum_{i=1}^{n-1} x^i \frac{\partial^2\phi}{\partial x_i\partial x_n}+\frac{\partial\phi}{\partial x_n} \ \ {\rm{on}}\ \ T_r,$$
so that 
\begin{equation}\label{I1-b}
\Big|I_{1b}\Big|\lesssim \int_{T_r}(|\nabla \phi|^2 |x||\nabla^2\phi|+|\nabla\phi|)\lesssim r^{n-1}.
\end{equation}
For $I_{1a}$, 
we have
\begin{eqnarray}
\Big|I_{1a}\Big|&=&\Big|\int_{S_r^+}\Big\langle{\partial\over\partial r}\Big(\mathbb B(u)(\nabla u,\nabla u)\Big),
\ x\cdot\nabla\phi\Big\rangle\Big|\nonumber\\
&\lesssim& r\int_{S_r^+}\Big(|\nabla u|^3+|\nabla u||\nabla^2u|\Big).\label{I1-a}
\end{eqnarray}
Putting the estimates for $I_{1a}$, $I_{1b}$, and $I_2$ together, we obtain

\begin{equation}
\Big|I\Big|
\lesssim \Big[\int_{B_r^+}\Big(|\nabla u|^3+|\nabla u||\nabla^2 u|\Big)+
r\int_{S_r^+}\Big(|\nabla u|^3+|\nabla u||\nabla^2u|\Big)+r^{n-1}\Big]. \label{I-est}
\end{equation}
For $III$, since 
$$\Delta(\mathbb P(u))=D^2\mathbb P(u)(\nabla u,\nabla u)+D\mathbb P(u)(\Delta u),$$
we have
\begin{equation}\label{I3-est}
\Big|III\Big|\lesssim \int_{B_r^+}\Big(|\nabla u|^4+|\Delta u|^2\Big)|x||\nabla \phi|
\lesssim r\int_{B_r^+}\Big(|\nabla u|^4+|\Delta u|^2\Big).
\end{equation}
For $II$,  we use integration by parts to estimate
\begin{eqnarray*}
II&=&2\int_{B_r^+}\nabla\cdot\Big\langle\langle\nabla(\mathbb P(u)),
\ \Delta u\rangle, \ x\cdot\nabla\phi\Big\rangle\\
&-&2\int_{B_r^+}\Big\langle\langle \nabla (\mathbb P(u)),\  \Delta u\rangle,\ \nabla(x\cdot\nabla\phi)\Big\rangle\\
&=&-2\int_{B_r^+}\Big\langle\langle\nabla(\mathbb P(u)),\ \Delta u\rangle,\ \nabla(x\cdot\nabla\phi)\Big\rangle\\
&&+2\int_{S_r^+}\Big\langle\Big\langle{\partial\over\partial r}(\mathbb P(u)),\ \Delta u\Big\rangle,\ x\cdot\nabla\phi\Big\rangle\\
&&-2\int_{T_r}\Big\langle\Big\langle\frac{\partial}{\partial x_n}(\mathbb P(u)),\ \Delta u\Big\rangle, \ x\cdot\nabla\phi\Big\rangle.
\end{eqnarray*}
Since 
$$\Delta u=(\Delta u)^T+\mathbb B(u)(\nabla u, \nabla u),$$
where $(\Delta u)^T=\mathbb P(u)(\Delta u)\in T_u N$ is the tangential component of $\Delta u$, 
and
$$\mathbb B(u)(\frac{\partial u}{\partial x_n},\ (\Delta u)^T)\perp x\cdot\nabla\phi \  {\rm{on}}\ T_r,$$
we have that, on $T_r$,
\begin{eqnarray*}
&&\Big\langle\Big\langle\frac{\partial}{\partial x_n}(\mathbb P(u)),\ \Delta u\Big\rangle, \ x\cdot\nabla\phi\Big\rangle\\
&=&-\Big\langle \mathbb B(u)(\frac{\partial u}{\partial x_n},\ (\Delta u)^T), \ x\cdot\nabla\phi\Big\rangle+\Big\langle\Big\langle D\mathbb P(u)(\frac{\partial u}{\partial x_n}),
\ \mathbb B(u)(\nabla u,\nabla u)\Big\rangle, \ x\cdot\nabla\phi\Big\rangle\\
&=&\Big\langle\Big\langle D\mathbb P(u)(\frac{\partial u}{\partial x_n}),
\ \mathbb B(u)(\nabla u,\nabla u)\Big\rangle, \ x\cdot\nabla\phi\Big\rangle\\
&=& \Big\langle\Big\langle D\mathbb P(\phi)(\frac{\partial \phi}{\partial x_n}),
\ \mathbb B(\phi)(\nabla \phi),\nabla \phi)\Big\rangle, \ x\cdot\nabla\phi\Big\rangle
\end{eqnarray*}
Therefore we have
\begin{eqnarray}
\Big|II\Big|&\lesssim&
\Big[\int_{B_r^+} |\nabla u||\Delta u|(|\nabla\phi|+|x||\nabla^2\phi|)
+\int_{S_r^+}|\frac{\partial u}{\partial r}||\Delta u||x||\nabla\phi|
+\int_{T_r}|x||\nabla\phi|^4\Big]\nonumber\\
&\lesssim&
\Big[\int_{B_r^+} |\nabla u||\Delta u|
+r\int_{S_r^+}|\nabla u||\Delta u|+r^n\Big] .\label{II-est}
\end{eqnarray}
Putting (\ref{I-est}), (\ref{I3-est}), and (\ref{II-est}) together and applying H\"older and Young's inequalities,
we obtain
\begin{eqnarray} \label{lastterm-est}
&&\Big|2r^{3-n}\int_{B_r^+}\Big\langle\Delta^2 u, \ x\cdot\nabla\phi\Big\rangle\Big|\nonumber\\
&\lesssim&
\Big[r^{3-n}\int_{B_r^+}\Big[r(|\nabla u|^4+|\nabla^2u|^2)+(|\nabla u|^3+|\nabla u||\nabla^2 u|)\Big]\nonumber\\
&+& r^{4-n}\int_{S_r^+}\left(|\nabla u|^3+|\nabla u||\nabla^2 u|\right)+r^3\Big]\nonumber\\
&\lesssim& 1+ r^{4-n}\int_{B_r^+}(|\nabla u|^4+|\nabla^2 u|^2)
+\frac{d}{dr}\left(r^{4-n}\int_{B_r^+}|\nabla u|^3+|\nabla u||\nabla^2 u|\right).
\end{eqnarray}

Similar to the derivation of interior monotonicity formula by Chang-Wang-Yang \cite{CWY}  page 1123-1124, and Angelsberg \cite{A} page 291-292, we estimate the first term of (\ref{radial3}) as follows.
\begin{eqnarray}
&&2r^{3-n}\int_{S_r^+}\Big[\left\langle\Delta (u-\phi),\ {\partial\over \partial r}\left(x\cdot\nabla (u-\phi)\right)\right\rangle-\left\langle{\partial\over \partial r}(\Delta (u-\phi)),\ x\cdot\nabla(u-\phi)\right\rangle\Big]\nonumber\\
&=&
2r^{2-n}\int_{S_r^+}\left[x^ix^j (u-\phi)_{ij}(u-\phi)_{kk} +x^i (u-\phi)_i (u-\phi)_{kk}
-x^ix^j (u-\phi)_j  (u-\phi)_{kki}\right]\nonumber\\
&=&IV(r).\nonumber
\end{eqnarray}
We estimate $IV$ term by term by integrating over $[\rho, r]$ as follows. For the first term of $IV$, 
applying integration by parts twice and using the fact that
$$(u-\phi)_n(u-\phi)_{ij}x^ix^j=0,
 \ (u-\phi)_k (u-\phi)_{ik} x^i x^n=0 \ \ {\rm{on}}\ T_1,$$
we obtain
\begin{eqnarray*}
\int_{B_r^+\setminus B_\rho^+} \frac{x^ix^j(u-\phi)_{ij}(u-\phi)_{kk}}{|x|^{n-2}}
&=&\int_{S_r^+\setminus S_\rho^+}\left(\frac{(u-\phi)_k (u-\phi)_{ij}x^ix^jx^k}{|x|^{n-1}}
-\frac{(u-\phi)_k (u-\phi)_{ik} x^i}{|x|^{n-3}}\right)\\
&+&\int_{B_r^+\setminus B_\rho^+}\left(\frac{(u-\phi)_i (u-\phi)_{ij}x^j}{|x|^{n-2}}+
\frac{(u-\phi)_{jk}(u-\phi)_{ik}x^i x^j}{|x|^{n-2}}\right)\\
&+&(n-2)\int_{B_r^+\setminus B_\rho^+}\frac{(u-\phi)_k(u-\phi)_{ij}x^ix^jx^k}{|x|^n}.
\end{eqnarray*}
For the second term of $IV$,  since $(u-\phi)_n x^i (u-\phi)_i=0$ on $T_1$, we have
\begin{eqnarray*}
\int_{B_r^+\setminus B_\rho^+}\frac{x^i (u-\phi)_i (u-\phi)_{kk}}{|x|^{n-2}}
&=&\int_{S_r^+\setminus S^+_\rho}\frac{|x^k(u-\phi)_k|^2}{|x|^{n-1}}
-\int_{B_r^+\setminus B_\rho^+}\frac{|\nabla (u-\phi)|^2}{|x|^{n-2}}\\
&+&\int_{B_r^+\setminus B_\rho^+}\left((n-2)\frac{|x^k (u-\phi)_k|^2}{|x|^n}
-\frac{x^i(u-\phi)_k (u-\phi)_{ik}}{|x|^{n-2}}\right).
\end{eqnarray*}
For the third term, applying integraton by parts and using the fact that
$$(u-\phi)_{ni} (u-\phi)_j x^i x^j=0,
\ \ {\rm{on}} \ \ T_1,$$
we have
\begin{eqnarray*}
&&\int_{B_r^+\setminus B_\rho^+} -\frac{x^i x^j (u-\phi)_j (u-\phi)_{kki}}{|x|^{n-2}}\\
&=&-\int_{S_r^+\setminus S_\rho^+}\frac{x^i x^j x^k (u-\phi)_{ki} (u-\phi)_j}{|x|^{n-1}}\\
&+&\int_{B_r^+\setminus B_\rho^+} \left(\frac{x^j(u-\phi)_k (u-\phi)_{kj}}{|x|^{n-2}}
+\frac{x^j (u-\phi)_j (u-\phi)_{kk}}{|x|^{n-2}}\right)\\
&+&\int_{B_r^+\setminus B_\rho^+}
\left((2-n)\frac{x^i x^j x^k (u-\phi)_j (u-\phi)_{ik}}{|x|^n}
+\frac{x^i x^j(u-\phi)_{ik}(u-\phi)_{jk}}{|x|^{n-2}}\right).
\end{eqnarray*}
Putting these identities together, we get
\begin{eqnarray}
&&\int_{\rho}^r IV(\tau)\,d\tau\nonumber\\
&=&2\int_{S_r^+\setminus S_\rho^+}\left(\frac{|x^k(u-\phi)_k|^2}{|x|^{n-1}}
-\frac{(u-\phi)_i (u-\phi)_{ij} x^j}{|x|^{n-3}}\right)\nonumber\\
&+& 2\int_{B_r^+\setminus B_\rho^+}
\Big(\frac{(u-\phi)_j (u-\phi)_{ij} x^i}{|x|^{n-2}}
+2\frac{|(u-\phi)_{ij}x^i|^2}{|x|^{n-2}}+\frac{(u-\phi)_i (u-\phi)_{jj} x^i}{|x|^{n-2}}\nonumber\\
&&-\frac{|\nabla (u-\phi)|^2}{|x|^{n-2}}+(n-2)\frac{|(u-\phi)_ix^i|^2}{|x|^n}\Big).
 \label{IV-identity1}
\end{eqnarray}
Using the identities
\begin{equation} 
2\int_{B_r^+\setminus B_\rho^+}\left(\frac{|\nabla (u-\phi)|^2}{|x|^{n-2}}
+\frac{(u-\phi)_j (u-\phi)_{ij} x^i}{|x|^{n-2}}\right)
=\int_{S_r^+\setminus S_\rho^+}\frac{|\nabla (u-\phi)|^2}{|x|^{n-3}},
\end{equation}
and
\begin{eqnarray}
&&\int_{B_r^+\setminus B_\rho^+}\left(\frac{(u-\phi)_j (u-\phi)_{ij} x^i
+|\nabla (u-\phi)|^2+(u-\phi)_i(u-\phi)_{jj}x^i}{|x|^{n-2}}
+(2-n)\frac{|x^i (u-\phi)_i|^2}{|x|^n}\right)\nonumber\\
&=&\int_{S_r^+\setminus S_\rho^+}\frac{|x^i (u-\phi)_i|^2}{|x|^{n-1}},
\end{eqnarray}
we then obtain
\begin{eqnarray}
\int_\rho^r IV(\tau)\,d\tau&=&2\int_{S_r^+\setminus S_\rho^+}\left(2\frac{|x^k(u-\phi)_k|^2}{|x|^{n-1}}-\frac{(u-\phi)_i (u-\phi)_{ij} x^j}{|x|^{n-3}}-2\frac{|\nabla (u-\phi)|^2}{|x|^{n-3}}\right)\nonumber\\
&+& 4\int_{B_r^+\setminus B_\rho^+}
\left(\frac{|(u-\phi)_j+x^i (u-\phi)_{ij}|^2}{|x|^{n-2}}+(n-2)\frac{|x^i(u-\phi)_i|^2}{|x|^n}\right). 
\label{IV-identity2}
\end{eqnarray}
Differentiating with respect to $r$, we obtain
\begin{eqnarray}
IV(r)&\ge& 2\frac{d}{dr}\int_{S_r^+}\left(2\frac{|x^k(u-\phi)_k|^2}{|x|^{n-1}}
-\frac{(u-\phi)_i (u-\phi)_{ij} x^j}{|x|^{n-3}}-2\frac{|\nabla(u-\phi)|^2}{|x|^{n-3}}\right)\nonumber\\
&+& 4\int_{S_r^+}\left(\frac{|(u-\phi)_j+x^i (u-\phi)_{ij}|^2}{|x|^{n-2}}+(n-2)\frac{|x^i(u-\phi)_i|^2}{|x|^n}\right)
 \label{IV-identity3}
\end{eqnarray}
Putting (\ref{thirdterm-est}),  (\ref{lastterm-est}), and  (\ref{IV-identity3})  into (\ref{radial3}), we obtain
\begin{eqnarray}\label{bdry_mono1}
&&{d\over dr}\left(r^{4-n}\int_{B_r^+}|\Delta (u-\phi)|^2\right)\nonumber\\
&\ge &-C\frac{d}{dr}\left(r^{4-n}\int_{B_r^+}|\nabla u|^3+|\nabla u||\nabla^2 u|\right)\nonumber\\
&+&2\frac{d}{dr}\int_{S_r^+}\left(2\frac{|x^k(u-\phi)_k|^2}{r^{n-1}}-\frac{(u-\phi)_i (u-\phi)_{ij} x^j}{r^{n-3}}
-2\frac{|\nabla (u-\phi)|^2}{r^{n-3}}\right)\nonumber\\
&+& 4\int_{S_r^+}\left(\frac{|(u-\phi)_j+x^i (u-\phi)_{ij}|^2}{|x|^{n-2}}+(n-2)\frac{|x^i(u-\phi)_i|^2}{|x|^n}\right)\nonumber\\
&-&C-Cr^{4-n}\int_{B_r^+}(|\nabla u|^4+|\nabla^2 u|^2).
\end{eqnarray}
It is easy to see
\begin{eqnarray}\label{l4-est1}
r^{4-n}\int_{B_r^+}|\nabla^2 u|^2
&\lesssim& r^{4-n}\int_{B_r^+}|\nabla^2 (u-\phi)|^2+ r^{4-n}\int_{B_r^+}|\nabla^2 \phi|^2\nonumber\\
&\lesssim& r^{4-n}\int_{B_r^+}|\nabla^2 (u-\phi)|^2+ r^4.
\end{eqnarray}\label{l4-est2}
By Nirenberg's interpolation inequality, we have
\begin{eqnarray}
r^{4-n}\int_{B_r^+}|\nabla u|^4
&\lesssim&
r^{4-n}\int_{B_r^+}|\nabla (u-\phi)|^4+r^{4-n}\int_{B_r^+}|\nabla \phi|^4\nonumber\\
&\lesssim& r^{4-n}\Big\|u-\phi\Big\|_{L^\infty(B_r^+)}^2 \int_{B_r^+}|\nabla^2(u-\phi)|^2+r^4\nonumber\\
&\lesssim& r^{4-n}\int_{B_r^+}|\nabla^2 (u-\phi)|^2+r^4. \label{hessian-est2}
\end{eqnarray}
Recall that the Bochner identity
$$\Delta|\nabla(u-\phi)|^2=2|\nabla^2(u-\phi)|^2+2\langle\nabla\Delta(u-\phi),\ \nabla(u-\phi)\rangle.$$
implies, after integrating over $B_r^+$ with integration by parts, that
\begin{eqnarray}
\frac{d}{dr}\left(r^{4-n}\int_{B_r^+}|\nabla^2(u-\phi)|^2\right)
&=&\frac{d}{dr}\left(r^{4-n}\int_{B_r^+}|\Delta(u-\phi)|^2\right)\nonumber\\
&+&2\frac{d}{dr}\left(r^{4-n}\int_{S_r^+}\left\langle\nabla(u-\phi), \frac{\partial}{\partial r}(\nabla(u-\phi))\right\rangle\right)\nonumber\\
&-&2\frac{d}{dr}\left(r^{4-n}\int_{S_r^+}\left\langle\Delta(u-\phi), \frac{\partial}{\partial r}(u-\phi)\right\rangle\right). \label{bochner}
\end{eqnarray}
Denote
$$f(r):=r^{4-n}\int_{S_r^+}\left\langle\nabla(u-\phi), \frac{\partial}{\partial r}(\nabla(u-\phi))\right\rangle
-r^{4-n}\int_{S_r^+}\left\langle\Delta(u-\phi), \frac{\partial}{\partial r}(u-\phi)\right\rangle,$$
$$g(r):=\int_{S_r^+}\left(2\frac{|x^k(u-\phi)_k|^2}{r^{n-1}}-\frac{(u-\phi)_i (u-\phi)_{ij} x^j}{r^{n-3}}
-2\frac{|\nabla (u-\phi)|^2}{r^{n-3}}\right).$$
Substituting (\ref{l4-est1}), (\ref{l4-est2}), and (\ref{bochner}) into (\ref{bdry_mono1}), we obtain
\begin{eqnarray} \label{bdry_mono2}
&&\frac{d}{dr}\Big[r^{4-n}\int_{B_r^+}|\nabla^2(u-\phi)|^2+C
\left(r^{4-n}\int_{B_r^+}|\nabla u|^3+|\nabla u||\nabla^2 u|\right)+Cr\Big]\nonumber\\
&&\ge 2 \frac{d}{dr} \Big(f+g\Big)-Cr^{4-n}\int_{B_r^+}|\nabla^2(u-\phi)|^2\nonumber\\
&+&4\int_{S_r^+}\left(\frac{|(u-\phi)_j+x^i (u-\phi)_{ij}|^2}{|x|^{n-2}}+(n-2)\frac{|x^i(u-\phi)_i|^2}{|x|^n}\right).
\end{eqnarray}
Notice that
\begin{eqnarray*} f(r)&\leq& C r^{4-n}\int_{S_r^+}|\nabla(u-\phi)||\nabla^2(u-\phi)|\\
&\le &C\Big[\frac{d}{dr}\Big(r^{4-n}\int_{B_r^+}|\nabla(u-\phi)||\nabla^2(u-\phi)|\Big)
+r^{4-n} \int_{B_r^+}|\nabla^2(u-\phi)|^2\Big]
\end{eqnarray*}
and
\begin{eqnarray*}
g(r)&\leq& C\Big[r^{3-n}\int_{S_r^+}|\nabla(u-\phi)|^2+r^{4-n}\int_{S_r^+}|\nabla(u-\phi)|\nabla^2(u-\phi)|\Big] \\
&\le& C\Big[\frac{d}{dr}\left(r^{3-n}\int_{B_r^+}|\nabla(u-\phi)|^2+r^{4-n}\int_{B_r^+}
|\nabla(u-\phi)||\nabla^2(u-\phi)|\right)\\
&&+r^{4-n}\int_{B_r^+}|\nabla^2(u-\phi)|^2\Big].
\end{eqnarray*}
Therefore, for sufficiently large $C>0$, we have
\begin{eqnarray}\label{bdry_mono3}
&&\frac{d}{dr}\Big[e^{Cr}r^{4-n}\int_{B_r^+}|\nabla^2(u-\phi)|^2+Ce^{Cr}
\left(r^{4-n}\int_{B_r^+}|\nabla u|^3+|\nabla u||\nabla^2 u|\right)+Cre^{Cr}\Big]\nonumber\\
&\ge& \frac{d}{dr}\Big[2e^{Cr}(f+g)\Big]
-\frac{d}{dr}\Big[Ce^{cr}\Big(r^{3-n}\int_{B_r^+}|\nabla(u-\phi)|^2+r^{4-n}\int_{B_r^+}
|\nabla(u-\phi)||\nabla^2(u-\phi)|\Big)\Big]\nonumber\\
&+&e^{Cr}\int_{S_r^+}\left(\frac{|(u-\phi)_j+x^i (u-\phi)_{ij}|^2}{|x|^{n-2}}+(n-2)\frac{|x^i(u-\phi)_i|^2}{|x|^n}\right).
\end{eqnarray}
Integrating (\ref{bdry_mono3}) for $r\in [\rho,R]$, we obtain
\begin{eqnarray}\label{bdry_mono4}
&&\rho^{4-n}\int_{B_\rho^+}|\nabla^2(u-\phi)|^2+\int_{B_R^+\setminus B_\rho^+}\left(\frac{|(u-\phi)_j+x^i (u-\phi)_{ij}|^2}{|x|^{n-2}}+(n-2)\frac{|x^i(u-\phi)_i|^2}{|x|^n}\right)
\nonumber\\
&\le& CRe^{CR}+e^{CR}R^{4-n}\int_{B_R^+}|\nabla^2(u-\phi)|^2
+Ce^{CR} R^{4-n}\int_{B_R^+}(|\nabla u|^3+|\nabla u||\nabla^2 u|)\nonumber\\
&+&Ce^{CR}\Big(R^{3-n}\int_{B_R^+}|\nabla(u-\phi)|^2+R^{4-n}\int_{B_R^+}|\nabla(u-\phi)||\nabla^2(u-\phi)|\Big)\nonumber\\
&+&Ce^{CR}\int_{S_R^+}\Big(R^{4-n}|\nabla (u-\phi)||\nabla^2(u-\phi)|+R^{3-n}|\nabla(u-\phi)|^2\Big)
\\
&+&Ce^{C\rho}\int_{S_\rho^+}\Big(\rho^{4-n}|\nabla (u-\phi)||\nabla^2(u-\phi)|+\rho^{3-n}|\nabla(u-\phi)|^2\Big)+e^{CR} g(R)-e^{cr}g(r)\nonumber.
\end{eqnarray}
Notice that by Poincar\'e inequlaity and H\"older inequality, we have
$$\Big(R^{3-n}\int_{B_R^+}|\nabla(u-\phi)|^2+R^{4-n}\int_{B_R^+}|\nabla(u-\phi)||\nabla^2(u-\phi)|\Big)
\leq C R^{5-n}\int_{B_R^+}|\nabla^2(u-\phi)|^2,$$
and
$$
R^{4-n}\int_{B_R^+}(|\nabla u|^3+|\nabla u||\nabla^2 u|)
\leq CR^{5-n}\int_{B_R^+}|\nabla^2(u-\phi)|^2+CR.
$$
Putting these two inequalities into (\ref{bdry_mono4}), we obtain
\begin{eqnarray}\label{bdry_mono50}
&&\rho^{4-n}\int_{B_\rho^+}|\nabla^2(u-\phi)|^2
+\int_{B_R^+\setminus B_\rho^+}\left(\frac{|(u-\phi)_j+x^i (u-\phi)_{ij}|^2}{|x|^{n-2}}+(n-2)\frac{|x^i(u-\phi)_i|^2}{|x|^n}\right)\nonumber\\
&\le& CRe^{CR}+(1+CR)e^{CR}R^{4-n}\int_{B_R^+}|\nabla^2(u-\phi)|^2\nonumber\\
&+&Ce^{CR}\Big[\left(R^{3-n}\int_{S_R^+}|\nabla (u-\phi)|^2\right)^\frac12
\left(R^{5-n}\int_{S_R^+}|\nabla^2(u-\phi)|^2\right)^\frac12\nonumber\\
&+&R^{3-n}\int_{S_R^+}|\nabla(u-\phi)|^2\Big]\nonumber\\
&+&Ce^{C\rho}\Big[\left(\rho^{3-n}\int_{S_\rho^+}|\nabla (u-\phi)|^2\right)^\frac12
\left(\rho^{5-n}\int_{S_\rho^+}|\nabla^2(u-\phi)|^2\right)^\frac12\nonumber\\
&+&\rho^{3-n}\int_{S_\rho^+}|\nabla(u-\phi)|^2\Big]+e^{CR}g(R)-e^{Cr}g(r).
\end{eqnarray}
It is clear that (\ref{bdry_mono50}) implies (\ref{bdry_mono}). \qed

\setcounter{section}{2} \setcounter{equation}{0}
\section{Boundary decay Lemma}

In this section, we will establish the bounday decay estimate for biharmonic maps that satisfy
the smallness condition (\ref{bdry_small}).

First we need to recall some notations.
For an open set $O\subset\mathbb R^n$, $1\le p<+\infty$ and
$0<\lambda\le n$, the Morrey space $M^{p,\lambda}(O)$ is defined by
$$
M^{p,\lambda}(O)=\Big\{
f\in L^p(O): \ \|f\|_{M^{p,\lambda}(O)}^p
=\sup_{B_r\subset O} \Big\{ r^{\lambda-n}\int_{B_r}|f|^p\Big\}<+\infty\Big\}.$$
Recall also that ${\rm{BMO}(O)}$ is defined by
$$
{\rm{BMO}}(O)=\Big\{
f\in L^1_{\rm{loc}}(O): \
\left[f\right]_{{\rm{BMO}}(O)}=\sup_{B_r\subset O}
\Big\{ r^{-n}\int_{B_r}|f-f_r|\Big\}<+\infty\Big\},$$
where $f_r=\frac{1}{|B_r|}\int_{B_r} f$ is the average of $f$ over $B_r$.

\begin{lemma} \label{bdry_decay_lem}
There exist $\epsilon_0>0$ and $\theta_0\in(0,\frac{1}{2})$ such that if
 $u\in W^{2,2}(B^{+}_1,N)$ is a biharmonic map satisfying
 \begin{equation}
 u\Big|_{T_1}=\phi\Big|_{T_1}
 \ {\rm{and}}\  \frac{\partial{u}}{\partial{x_n}}\Big|_{T_1}=\frac{\partial{\phi}}{\partial{x_n}}\Big|_{T_1}
\  \ {\rm{for\ some}}\ \phi\in C^{\infty}\Big(\overline{B^+_1},N\Big), \label{bdry_cond1}
\end{equation}
and
 \begin{equation}
\left\|\nabla^2u\right\|_{M^{2,4}(B^+_1)}+\left\|\nabla u\right\|_{M^{4,4}(B^+_1)}\leq \epsilon_0. \label{bdry_small}
 \end{equation}
 Then
  \begin{equation} \label{bdry_decay}
\Big\|\nabla u\Big\|_{M^{2,2}(B^+_{\theta_0})}\leq \frac{1}{2}
\Big\|\nabla u\Big\|_{M^{2,2}(B^+_1)}+C\Big\|\nabla\phi\Big\|_{C^1(B^+_1)}\theta_0 .
\end{equation}
In particular, $u\in C^\infty\Big(\overline{B_\frac12^+}, N\Big).$
\end{lemma}

\pf 
A crucial step  to establish (\ref{bdry_decay}) is the following claim.\\

\noindent{\bf Claim 1}.  For any $0<r\le\frac12$ and $\theta\in (0,\frac14)$, it holds
\begin{eqnarray}\label{bdry_decay11}
\Big\{(\theta r)^{2-n}\int_{B_{\theta r}^+}|\nabla u|^2\Big\}^\frac12
&\le& C\Big[\theta r (\|\nabla\phi\|_{C^1(B_{2r}^+)}+\|\nabla^2 u\|_{M^{2,2}(B_{2r}^+)})
+\theta^{1-\frac{n}2}\epsilon_0\|\nabla u\|_{M^{2,2}(B_{2r}^+)}\nonumber\\
&&+\|\nabla u\|_{M^{4,4}(B_{2r}^+)}^2\Big].
\end{eqnarray}

By scalings, it suffices to prove claim 1 for $r=\frac12$.  We will divide the proof into several steps.

\noindent{\it Step} 1.  First, we follow Struwe's scheme (see \cite {SM} pages 250-251) 
to rewrite the biharmonic map equation (\ref{biharm1}) into:
\begin{equation}\label{biharm3}
 \Delta^2u=\Delta(D\cdot\nabla u)+{\rm{div}}(E\cdot\nabla u)+F\cdot\nabla u  \ \ \ \text{in}\ \ B^+_1 , 
\end{equation}
and
$$
F=G+\Delta\Omega ,\ \ \Omega=(\Omega^{ij}=w^idw^j-w^jdw^i)\in H^1(B^+_1,so(l)\otimes\wedge^1\mathbb{R}^l),$$
where $w=w(u)=\nu\circ u$, with $\nu$ a unit normal vector of $N$, and
the coefficient functions $D,E,G$ and $\Omega$ depend on $u$ and satisfy
\begin{eqnarray}
     \left\{\begin{array}{llllll}
     |D|+|\Omega|\les|\nabla u|,\\
     |E|+|\nabla D|+|\nabla \Omega|\les |\nabla^2u|+|\nabla u|^2,\\
     |G|\les

 |\nabla^2u||\nabla u|+|\nabla u|^3.
 \end{array}\right.
\end{eqnarray}

\noindent{\it Step} 2. Extend $u$ from $B^+_1$ to $B_1$ , denote as $\tilde{u}$ , such that
\begin{equation}
\|\nabla^2\widetilde{u}\|_{M^{2,4}(B_1)}+\|\nabla \widetilde{u}\|_{M^{4,4}(B_1)}
\lesssim\|\nabla^2u\|_{M^{2,4}(B^+_1)}+\|\nabla u\|_{M^{4,4}(B^+_1)}
\leq C\epsilon_0. \label{small_cond1} \end{equation}
Let $\widetilde{\Omega}\in H^1(B_1, so(l)\otimes\wedge^1(\mathbb R^l))$ be the extension of 
$\Omega$ to $B_1$ given by 
$$\widetilde{\Omega}^{ij}=w^i(\widetilde{u})d(w^j(\widetilde{u}))
-w^j(\widetilde{u})d(w^i(\widetilde{u})), \  1\le i, j\le l.$$
Then one can check that $\widetilde{\Omega}$ satisfies
\begin{equation}\label{small_omega}
 \left\|\nabla\widetilde{\Omega}\right\|_{M^{2,4}(B_1)}
+\left\|\widetilde{\Omega}\right\|_{M^{4,4}(B_1)}
+\left\|\widetilde\Omega\right\|_{M^{2,2}(B_1)}\leq C \epsilon_0.
\end{equation}
Now we apply the Coulomb gauge construction in Morrey spaces, see \cite{RT} or \cite{SM},
to obtain that there exist $P\in H^2(B_1,SO(l))$
and $\xi\in H^2(B_1,so(l)\otimes\wedge^{n-2}\mathbb{R}^l)$ such that
\begin{eqnarray}
\left\{
\begin{array}{lllll}
d PP^{-1}+P\widetilde{\Omega}P^{-1}=*d\xi \ \ \text{in} \ \ B_1,\\
d(*\xi)=0 \ \ \text{in}\ \  B_1, \ \xi\Big|_{\partial B_1}=0.
\end{array}
\right.
\end{eqnarray}
Furthermore, $P$ and $\xi$ satisfy
\begin{equation}\label{gauge_morrey_est1}
\Big\|\nabla P\Big\|_{M^{2,2}(B_1)}\lesssim \left\|\widetilde\Omega\right\|_{M^{2,2}(B_1)}\lesssim
\Big\|\nabla u\Big\|_{M^{2,2}(B_1^+)}, \end{equation}
\begin{equation}\label{gauge_morrey_est2}
\ \Big\|\nabla P\Big\|_{M^{4,4}(B_1)}\lesssim \left\|\widetilde\Omega\right\|_{M^{4,4}(B_1)}
\lesssim \Big\|\nabla u\Big\|_{M^{4,4}(B_1^+)},
\end{equation}
and
\begin{eqnarray}\label{gauge_morrey_est3}
&&\left\|\nabla^2P\right\|_{M^{2,4}(B_1)}+
\left\|\nabla^2\xi\right\|_{M^{2,4}(B_1)}+\left\|\nabla\xi\right\|_{M^{4,4}(B_1)}
\nonumber\\&&
\lesssim \left\|\nabla\widetilde{\Omega}\right\|_{M^{2,4}(B_1)}
+\left\|\widetilde{\Omega}\right\|_{M^{4,4}(B_1)}\lesssim\epsilon_0.
\end{eqnarray}

\noindent{\it Step} 3.  Apply $P$ to the equation (\ref{biharm3}) (see, e.g. \cite{SM} page 254),
 we obtain
\begin{equation}
\Delta(P\Delta u)={\rm{div}}^2(D_P\otimes\nabla u)+{\rm{div}}(E_P\cdot\nabla u)+G_P\cdot\nabla u
+*d\Delta\xi\cdot P\nabla u\ \ {\rm{on}}\ B_1^+, \label{biharm4}
\end{equation}
where $D_P, E_P$, and $G_P$ satisfy:
\begin{eqnarray}
\left\{
\begin{array}{lllll}
|D_P|\lesssim(|\nabla u|+|\nabla P|),\\
|\nabla D_P|+|E_P|\lesssim(|\nabla^2u|+|\nabla u|^2+|\nabla^2P|+|\nabla P|^2),\\
|G_P|\lesssim(|\nabla^2u|+|\nabla^2P|)(|\nabla u|+|\nabla P|)+(|\nabla u|^3+|\nabla P|^3).
\end{array}
\right.
\end{eqnarray}

Now we need 

\noindent{\bf Claim} 2.  For sufficiently small $\epsilon_0>0$,  there exists $P_{\frac34}\in SO(l)$ such
that for any $1\le q<\infty$, there exists $C_q>0$ so that
\begin{equation}\label{const_p}
\left(\frac{1}{|B_{\frac34}|}\int_{B\frac34}|P-P_{\frac34}|^q\right)^{\frac{1}{q}}
\leq C_q\left\|\nabla P\right\|_{M^{1,1}(B_1)}.
\end{equation}

To show (\ref{const_p}), first observe that by H\"older's inequality, we have
$$\left\|\nabla P\right\|_{M^{1,1}(B_1)}\leq\left\|\nabla P\right\|_{M^{4,4}(B_1)}\leq\epsilon_0.$$
By Poincar\'e's inequality, we have that
$P\in {\rm{BMO}}(B_1)$ and 
$$\left[P\right]_{{\rm{BMO}}(B_1)}\lesssim\epsilon_0.$$
Hence, by John-Nirenberg's inequality, we conclude that for any $1\le q<\infty$, there is $C_q>0$ such
that
\begin{equation}\label{bmo-small1}
\left(\frac{1}{|B_{\frac34}|}\int_{B_{\frac34}}|P-(P)_{B_{\frac34}}|^q\right)^{\frac{1}{q}}
\leq C_q \left[P\right]_{{\rm{BMO}}(B_\frac34)}\leq C_q\left\|\nabla P\right\|_{M^{1,1}(B_1)},
\end{equation}
where $(P)_{B_{\frac34}}=\frac{1}{|B_{\frac34}|}\int_{B_{\frac34}}P$ is 
the average of $P$ over $B_{\frac34}$.

Now we need to project $(P)_{B_{\frac34}}$ into $SO(l)$. Since
\begin{equation*}
{\rm{dist}}((P)_{\frac34},\ SO(l))\leq\frac{1}{|B_{\frac{3}{4}}|}\int_{B_{\frac{3}{4}}}
|P-(P)_{B_{\frac{3}{4}}}|
\leq[P]_{{\rm{BMO}}(B_{\frac34})}\leq\epsilon_0,
\end{equation*}
there exists a $P_{\frac{3}{4}}\in SO(l)$ such that 
\begin{equation}
{\rm{dist}}((P)_{B_{\frac{3}{4}}},\ SO(l))=\left|P_{\frac{3}{4}}-(P)_{B_{\frac{3}{4}}}\right|
\leq \left[P\right]_{{\rm{BMO}}(B_\frac34)}\le C\left\|\nabla P\right\|_{M^{1,1}(B_1)}.\label{bmo-small2}
\end{equation}
Combining (\ref{bmo-small1}) with (\ref{bmo-small2}), we obtain (\ref{const_p}). 

By Fubini's theorem, we may assume that
$u\Big|_{\partial B^+_{\frac{3}{4}}}\in H^2(\partial B_{\frac34}^+)$.
Set $\omega:B^+_{\frac34}\rightarrow\mathbb{R}^l$ to be  the biharmonic extension of
$u$ on $\partial B_{\frac34}^+$, i.e., 
\begin{eqnarray}
\left\{
\begin{array}{lllll}
\triangle^2\omega=0\ \ \ \ \ \ \text{in}\ \ B^+_{\frac{3}{4}}\\
\omega=u\ \ \ \ \ \ \ \ \ \ \text{on}\ \  \partial B^+_{\frac{3}{4}}\\
\frac{\partial\omega}{\partial\nu}=\frac{\partial u}{\partial\nu}\ \ \ \ \ \ \ \ \text{on}
\ \ \partial B^+_{\frac{3}{4}}.
\end{array}
\right.
\end{eqnarray}
Since $(u,\nabla u)\Big|_{T_1}=(\phi, \nabla\phi)\in C^\infty(T_1)$,
it follows from both interior and boundary regularity theory (see, Agmon-Douglas-Nirenberg \cite{ADN})
that $\omega\in C^{\infty}\Big(\overline{B^+_{\frac{5}{8}}},\mathbb{R}^l\Big)$ and for $m\ge 1$
\begin{equation}\label{est-biharm-fun}
\left\|\omega\right\|_{C^m\left(\overline{B^+_{\frac{5}{8}}}\right)}\lesssim \Big[\|\phi\|_{C^{m+1}\left(\overline{B^+_{\frac{3}{4}}}\right)}+
\left\|\nabla^2\omega\right\|_{L^2(B^+_{\frac{3}{4}})}\Big]
\lesssim\Big[\|\phi\|_{C^{m+1}\left(\overline{B^+_{\frac{3}{4}}}\right)}+
\left\|\nabla^2u\right\|_{L^2(B^+_{\frac{3}{4}})}\Big],
\end{equation}
since $\omega$ minimizes the Hessian energy, i.e.,
\begin{equation}
\int_{B^+_{\frac{3}{4}}}|\nabla^2\omega|^2\leq \int_{B^+_{\frac{3}{4}}}|\nabla^2 u|^2.
\end{equation}

Let $U$ be a bounded, smooth domain such that $B_\frac12^+\subset U\subset B_\frac34^+$. In particular,
$T_\frac12\subset\partial U$ and $\partial U\cap\partial\mathbb R^n_+\subset T_\frac34$.
To proceed with the proof, we recall the estimate of Green functions of 
$\triangle^2$ on $U$.
Let $G(x,y)$ be the Green function of $\triangle^2$ on $U$:
for $x\in U$, it holds
\begin{eqnarray}
\left\{
\begin{array}{lllll}
\triangle^2_y G(x,y)=\delta_x(y) \ \ \ \ \ \text{in}\ U\\
\ G(x,y)=0 \ \ \ \ \ \ \ \ \ \ \ \ \ \text{on}\ \ \partial U\\
\frac{\partial G}{\partial\nu_y}(x,y)=0\ \ \ \ \ \ \ \ \  \ \ \ \text{on}\ \ \partial U.
\end{array}
\right.
\end{eqnarray}
Then we have (see Dall'Acqua-Sweers \cite{DS} Theorem 3 and Theorem B.2, 
due to Krasovskii \cite{Kra}):
\begin{equation}\label{green-est}
\Big|D_x^{\alpha}D^{\beta}_yG(x,y)\Big|\lesssim\Big|x-y\Big|^{4-n-|\alpha|-|\beta|},\ \ x,
\ y\in U, \ \
\text{for any}\ \ |\alpha|,\ |\beta|\geq 0 .
\end{equation}

For $P_{\frac34}\in SO(l)$ given by Claim 1, we then have
\begin{eqnarray*}
\triangle(P\triangle u)
=\triangle(P\triangle u)-P_{\frac{3}{4}}\triangle^2\omega
=\triangle^2(P_{\frac{3}{4}}(u-\omega))+\triangle((P-P_{\frac{3}{4}})\triangle u) .
\end{eqnarray*}
For simplicity, we may assume $P_{\frac{3}{4}}=\mathbb I_l\in SO(l)$, the identity $l\times l$-matrix.
Then (\ref{biharm4}) gives
\begin{eqnarray}
\left\{
\begin{array}{lllll}
\triangle^2(u-\omega)=\triangle((\mathbb I_l-P)\triangle u)+{\rm{div}}^2(D_P\otimes\nabla u)+ {\rm{div}}(E_P\cdot\nabla u)\\
\ \ \ \ \ \ \  \ \ \ \ \ \ \ \ \ \ +G_P\cdot\nabla u+*d\triangle\xi\cdot P\nabla u \ \ \ \ \ \ \ \ \ \ \ \ \ \ \ \
\ \ \ \ \ \ \ \ \ \ \ \ \ \ \ \ \text{in}\ \ U\\
u-\omega=\frac{\partial}{\partial \nu}(u-\omega)=0\ \ \ \ \ \ \ \ \
 \ \ \ \ \ \ \ \ \ \ \ \ \ \ \ \ \ \ \ \ \ \ \ \ \ \ \ \ \ \ \ \ \ \ \ \ \ \ 
\ \ \ \ \text{on}\ \ \partial U.
\end{array}
\right.
\end{eqnarray}
Therefore, by the represenation formula 
\begin{equation}
(u-\omega)(x)=\int_{U}G(x,y)\triangle^2(u-\omega)(y)\,dy,
\end{equation}
we have that for any $x\in U$,
\begin{align}\label{represent}
(u-\omega)(x)&=\int_{U}\triangle_yG(x,y)((\mathbb I_l-P)\triangle u)(y)\\&
+\int_{U}[\nabla^2_yG(x,y)(D_P\otimes\nabla u)(y)-\nabla_yG(x,y)(E_P\cdot\nabla u)(y)]\nonumber\\&
+\int_{U}G(x,y)(G_P\cdot\nabla u)(y)-\int_{U}\triangle
\xi(y)d_y(G(x,y)P(y))\wedge du(y) \nonumber.
\end{align}
Differentiating (\ref{represent}) with respect to $x$, we have 
\begin{align}
|\nabla(u-\omega)(x)|&\lesssim\int_{U}|\nabla_x\nabla_y^2G(x,y)||\mathbb I_l-P||\triangle u|(y)\\&
+\int_{U}|\nabla_x\nabla_y^2G(x,y)||D_P||\nabla u|(y)\nonumber\\&
+\int_{U}|\nabla_x\nabla_yG(x,y)||E_P||\nabla u|(y)\nonumber\\&
+\int_{U}|\nabla_xG(x,y)|G_P||\nabla u|(y)\nonumber\\&
+\int_{U}|\nabla_xG(x,y)||\triangle\xi||\nabla P||\nabla u|(y)\nonumber\\&
+\int_{U}|\nabla_x\nabla_yG(x,y)||\Delta\xi||\nabla u|(y)\nonumber\\&
=A_1+A_2+A_3+A_4+A_5+A_6\nonumber
\end{align}
We now estimate $A_1,\cdots, A_6$ as follows. First recall the Riesz potential of order $\alpha$
$$I_{\alpha}(f)(x)=\int_{\mathbb{R}^n}|x-y|^{\alpha-n}f(y)\,dy, \ f\in L^1(\mathbb R^n),$$
for $0<\alpha<n$. Then we have
\begin{equation*}
\Big|A_1\Big|\lesssim I_1\Big(|\mathbb I_l-P||\triangle u|\chi_{B^+_{\frac34}}\Big),
\end{equation*}
where $\chi_{B_{\frac34}^+}$ is the characteristic function of $B_{\frac34}^+$.
By the standard estimate on Riesz potentials, we have
\begin{equation*}
\left\|A_1\right\|_{L^2(B^+_{\frac{1}{2}})}\lesssim\left\|\mathbb I_l-P\right\|_{L^n(B^+_{\frac34})}
\left\|\nabla^2u\right\|_{L^2(B^+_{\frac34})},
\end{equation*}
so that
\begin{align}
\left\|A_1\right\|_{L^2(B^+_{\frac{1}{2}})}&
\lesssim\left\|\nabla^2 u\right\|_{M^{2,4}(B_1^+)}\left\|\mathbb I_l-P\right\|_{L^n(B^+_{\frac34})}\\
&\lesssim\epsilon_0\|\nabla P\|_{M^{2,2}(B_1^+)}\lesssim\epsilon_0\Big\|\nabla u\Big\|_{M^{2,2}(B_1^+)}.\nonumber
\end{align}
For $A_2$, we have
\begin{align*}
|A_2|
\lesssim I_1\Big(|D_P||\nabla u|\chi_{B^+_{\frac34}}\Big)
\lesssim I_1\Big((|\nabla u|^2+|\nabla P|^2)\chi_{B^+_{\frac34}}\Big).
\end{align*}
Since $|\nabla u|^2+|\nabla P|^2\in M^{2,4}(B^+_1)$, one can check
\begin{equation*}
(|\nabla u|^2+|\nabla P|^2)\chi_{B^+_{\frac34}}\in M^{2,4}(\mathbb{R}^n).
\end{equation*}
Applying Adams' Morrey space estimate of Riesz potentials (see \cite{W1} Proposition 4.2), we obtain
\begin{equation*}
|I_1(|D_P||\nabla u|\chi_{B_\frac34^+})|\in M^{4,4}(\mathbb{R}^n),
\end{equation*}
and
\begin{align}
\|A_2\|_{M^{4,4}(B^+_{\frac{1}{2}})}&\le \|A_2\|_{M^{4,4}(\mathbb R^n)}\lesssim
\Big\||D_P||\nabla u|\chi_{B_\frac34^+})|\Big\|_{M^{2,4}(\mathbb{R}^n)}\nonumber\\
&\lesssim\Big [\|\nabla u\|_{M^{4,4}(B^+_{\frac34})}^2+\|\nabla P\|_{M^{4,4}(B^+_{\frac34})}^2\Big]
\lesssim \|\nabla u\|_{M^{4,4}(B^+_1)}^2.
\end{align}
For $A_3$, we have
\begin{align*}
|A_3|&
\lesssim I_2\Big(|E_P||\nabla u|\chi_{B^+_{\frac34}}\Big)
\lesssim I_2\Big((|\nabla^2 u|+|\nabla u|^2+|\nabla^2P|+|\nabla P|^2)|\nabla u|\chi_{B^+_\frac34}\Big).
\end{align*}
By H\"older's inequality, we have
\begin{equation*}
(|\nabla^2u|+|\nabla u|^2+|\nabla^2P|+|\nabla P|^2)|\nabla u|
\chi_{B^+_{\frac34}}\in M^{1,3}(\mathbb{R}^n)
\end{equation*}
Applying the weak estimate (4.6) in \cite{W1} page 430, we have
$A_3\in M^{3,3}_*(\mathbb{R}^n)$, and
\begin{align}
\Big\|A_3\Big\|_{M^{2,2}(B_\frac12^+)}&\leq\Big\|A_3\Big\|_{M^{3,3}_*(B^+_{\frac{1}{2}})}
\leq\Big\|A_3\Big\|_{M^{\frac43,4}(\mathbb{R}^n)}\nonumber\\
&\leq\Big\||\nabla^2 u|+|\nabla u|^2+|\nabla^2 P|+|\nabla P|^2\Big\|_{M^{2,4}(B^+_1)}
\Big\|\nabla u\Big\|_{M^{2,2}(B^+_1)}\nonumber\\&
\lesssim\epsilon_0\Big\|\nabla u\Big\|_{M^{2,2}(B^+_1)}
\end{align}
Here and below we use the fact that 
$\|f\|_{M^{2,2}(B^+_{\frac{1}{2}})}\leq\|f\|_{M_*^{p,p}(B^+_{\frac{1}{2}})}$ for any $2<p<+\infty$.
For $A_4$, we have
\begin{align*}
|A_4|&
\lesssim I_3\Big(|G_P||\nabla u|\chi_{B^+_{\frac{3}{4}}}\Big)\\&
\lesssim I_3\Big([(|\nabla^2 u|+|\nabla^2P|)(|\nabla u|+|\nabla P|)|\nabla u|
+(|\nabla u|^4+\nabla P|^3|\nabla u|)]\chi_{B^+_{\frac{3}{4}}}\Big).
\end{align*}
Observe 
\begin{equation*}
\Big[(|\nabla^2u|+|\nabla^2P|)(|\nabla u|+|\nabla P|)|\nabla u|
+(|\nabla u|^4+|\nabla P|^3|\nabla u|)\Big]\chi_{B^+_{\frac{3}{4}}}\in M^{1,4}(\mathbb{R}^n).
\end{equation*}
Applying the weak estimate (4.6) in \cite{W1} page 430, we get $A_4\in M_*^{4,4}(\mathbb R^n)$
and 
\begin{align*}
\Big\|A_4\Big\|_{M_*^{4,4}(B^+_{\frac{1}{2}})}&
\lesssim\Big [\|\nabla^2u\|_{M^{2,2}(B^+_{\frac{3}{4}})}+\|\nabla^2P\|_{M^{2,2}(B^+_{\frac{3}{4}})}\Big]
\cdot\Big [\|\nabla u\|^2_{M^{4,4}(B^+_{\frac{3}{4}})}+\|\nabla P\|^2_{M^{4,4}(B^+_{\frac{3}{4}})}\Big]\\&
+\Big[\|\nabla u\|^4_{M^{4,4}(B^+_{\frac{3}{4}})}+\|\nabla P\|^4_{M^{4,4}(B^+_{\frac{3}{4}})}\Big]\\&
\lesssim \|\nabla u\|_{M^{4,4}(B^+_1)}^2+
\|\nabla P\|_{M^{4,4}(B^+_1)}^2
\lesssim \|\nabla u\|_{M^{4,4}(B^+_1)}^2.
\end{align*}
Hence we obtain
\begin{equation}
\Big\|A_4\Big\|_{M^{2,2}(B^+_{\frac{1}{2}})}\lesssim
\|\nabla u\|_{M^{4,4}(B^+_1)}^2.
\end{equation}
For $A_5$ , since we have 
\begin{equation*}
|A_5|\lesssim I_3\Big(|\triangle\xi||\nabla P||\nabla u|\chi_{B^+_{\frac{3}{4}}}\Big),
\end{equation*}
similar to the estimate of $A_4$, we obtain
\begin{align*}
\Big\|A_5\Big\|_{M_*^{4,4}(B^+_{\frac{1}{2}})}&
\lesssim \|\triangle\xi\|_{M^{2,4}(B^+_{\frac{3}{4}})}\|\nabla P\|_{M^{4,4}(B^+_{\frac{3}{4}})}
\|\nabla u\|_{M^{4,4}(B^+_{\frac{3}{4}})}
\lesssim\Big\|\nabla u\Big\|_{M^{4,4}(B^+_1)}^2,
\end{align*}
so that
\begin{equation}
\Big\|A_5\Big\|_{M^{2,2}(B^+_{\frac{1}{2}})}\lesssim
\Big\|\nabla u\Big\|_{M^{4,4}(B^+_1)}^2.
\end{equation}
For $A_6$, we have
\begin{equation*}
|A_6|(x)\lesssim I_2\Big(|\triangle\xi||\nabla u|\chi_{B_\frac34^+}\Big).
\end{equation*}
Similar to the estimate of $A_3$,  since $|\triangle\xi||\nabla u|\chi_{B_\frac34^+}
\in M^{1,3}(\mathbb R^n)$, we get 
\begin{equation}
\Big\|A_6\Big\|_{M^{2,2}(B^+_{\frac{1}{2}})}
\le\Big\|A_6\Big\|_{M^{3,3}_*(\mathbb R^n)}
\lesssim\epsilon_0
\Big\|\nabla u\Big\|_{M^{2,2}(B^+_1)}.
\end{equation}
Putting all these estimates together, we then obtain
\begin{equation}\label{bdry_decay20}
\Big(\int_{B^+_\frac12}|\nabla(u-\omega)|^2\Big)^{\frac{1}{2}}
\lesssim\epsilon_0\Big\|\nabla u\Big\|_{M^{2,2}(B^+_1)}+
\Big\|\nabla u\Big\|_{M^{4,4}(B^+_1)}^2.
\end{equation}
On the other hand, by the estimate of biharmonic function (\ref{est-biharm-fun}) we have that
for any $0<\theta<\frac{1}{4}$,
\begin{equation}\label{est-biharm-fun1}
\theta^{2-n}\int_{B^+_{\theta}}|\nabla\omega|^2\lesssim \theta^2\Big(\|\nabla\phi\|_{C^1(B^+_\frac12)}^2+
\int_{B_\frac12^+}|\nabla^2 u|^2\Big).
\end{equation}
Combining (\ref{bdry_decay20}) with (\ref{est-biharm-fun1}),  we have 
\begin{align}\label{bdry_decay30}
\Big(\theta^{2-n}\int_{B^+_{\theta }}|\nabla u|^2\Big)^{\frac12}
&\lesssim \theta^{1-\frac{n}{2}}\Big[\epsilon_0\|\nabla u\|_{M^{2,2}(B^+_1)}+
\|\nabla u\|_{M^{4,4}(B^+_1)}^2\Big]\nonumber\\
&+\theta\Big(\|\nabla\phi\|_{C^1(B_1^+)}+\|\nabla^2 u\|_{M^{2,4}(B_1^+)}\Big).
\end{align}
It is clear that (\ref{bdry_decay11}) follows from (\ref{bdry_decay30}).

Next we indicate how to prove (\ref{bdry_decay}) by applying (\ref{bdry_decay11}) and the interior regularity
theorem (see \cite{W1}, \cite{SM}, \cite{CWY}).  This is summarized as the third claim.\\

\noindent{\bf Claim} 3. For $\theta\in (0,\frac14)$, it holds
\begin{align}\label{bdry_decay40}
\|\nabla u\|_{M^{2,2}(B^+_{\theta})}&
\lesssim \theta^{1-\frac{n}{2}}\left[\epsilon_0\|\nabla u\|_{M^{2,2}(B^+_1)}
+\|\nabla u\|_{M^{4,4}(B^+_1)}^2\right]\\&
+\theta\Big(\|\nabla\phi\|_{C^1(B_1^+)}+\|\nabla^2 u\|_{M^{2,4}(B_1^+)}\Big).\nonumber
\end{align}

To prove (\ref{bdry_decay40}), let $y=(y',y^n)\in B^+_{\theta}$ and $\tau>0$ such that
$B_{\tau}(y)\subset B_{\theta}$. We want to show (\ref{bdry_decay40}) holds
with left hand side replaced by
$\Big(\tau^{2-n}\int_{B_\tau(y)}|\nabla u|^2\Big)^\frac12$.  We divide it into three cases:
\\

\noindent (i) $y^n>0$ and $B_{\tau}(y)\subset B^+_{\theta}$ . In this case, 
since we have 
$$\|\nabla u\|_{M^{4,4}(B_\tau(y))}+\|\nabla^2 u\|_{M^{2,4}(B_\tau(y))}\le
\|\nabla u\|_{M^{4,4}(B_1^+)}+\|\nabla^2 u\|_{M^{2,4}(B_1^+)}\le\epsilon_0,$$
it follows from the interior regularity theorem for stationary biharmonic maps and (\ref{bdry_decay11}) that
for any $0<\alpha<1$,
\begin{align*}
\left(\tau^{2-n}\int_{B_{\tau}(y)}|\nabla u|^2\right)^{\frac{1}2}&\lesssim\left(\frac{\tau}{y^n}\right)^{\alpha}
\left((y^n)^{2-n}\int_{B_{y^n}(y)}|\nabla u|^2\right)^{\frac{1}{2}}\\&
\lesssim\left((2y^n)^{2-n}\int_{B_{2y^n}(y',0)}|\nabla u|^2\right)^{\frac{1}{2}}\\&
\lesssim \theta^{1-\frac{n}{2}}\Big[\epsilon_0\|\nabla u\|_{M^{2,2}(B^+_1)}
+\|\nabla u\|^2_{M^{4,4}(B^+_1)}\Big]+C\theta.
\end{align*} 

\noindent (ii) $y^n>0$ and $B_{\tau}(y)\setminus B^+_{\theta}\neq\emptyset$.
Then $y^n<\tau$
so that $B_{\tau}(y)\subset B_{2\tau}(y',0)$ (as $2\tau\leq2\theta<\frac{1}{2}$). Therefore,
by (\ref{bdry_decay}) we have
\begin{align*}
\Big(\tau^{2-n}\int_{B_{\tau}(y)\cap B^+_{\theta}}|\nabla u|^2)^{\frac{1}{2}}&
\lesssim\Big((2\tau)^{2-n}\int_{B_{2\tau}(y',0)}|\nabla u|^2\Big)^{\frac{1}{2}}\\&
\lesssim \theta^{1-\frac{n}{2}}\Big[\epsilon_0\|\nabla u\|_{M^{2,2}(B^+_1)}
+\|\nabla u\|^2_{M^{4,4}(B^+_1)}\Big]+C\theta.
\end{align*}

\noindent (iii) $y^n=0$.  By translation,  one can easily see
that (\ref{bdry_decay}) holds for balls with center $y=(y',0)$.\\

Now taking supremum over all such $y$ and $\tau$, we obtain (\ref{bdry_decay40}).
Finally, recall the following interpolation inequality (see \cite{SM} Proposition 3.2 or \cite{W1} Proposition 4.3):
\begin{align}\label{inter_ineq}
\|\nabla u\|^2_{M^{4,4}(B^+_1)}&\lesssim\|\nabla u\|_{M^{2,2}(B^+_1)}(\|\nabla^2u\|_{M^{2,4}(B^+_1)}
+\|\nabla u\|_{M^{4,4}(B^+_1)})\\
&\lesssim\epsilon_0\|\nabla u\|_{M^{2,2}(B^+_1)}.\nonumber
\end{align}
Substituting (\ref{inter_ineq})  into (\ref{bdry_decay40}), we obtain

\begin{equation}\label{bdry_decay50}
\Big\|\nabla u\Big\|_{M^{2,2}(B^+_{\theta})}\lesssim \theta^{1-\frac{n}{2}}
\epsilon_0\Big\|\nabla u\Big\|_{M^{2,2}(B^+_1)}+C\theta.
\end{equation}
By chosing $\epsilon_0>0$ and $\theta_0=\theta(\epsilon_0)$ sufficiently small, this yields
(\ref{bdry_decay}).

It is clear that repeated iterations of (\ref{bdry_decay40}) imply
that there exists $\alpha\in (0,1)$ such that
$$\Big\|\nabla u\Big\|_{M^{2,2}(B_r^+)}\le Cr^\alpha, \ \forall 0<r\le\frac12.$$
This, with Morrey's decay Lemma and the interior regularity theorem, yields $u\in C^\alpha\Big(\overline{B_\frac12^+}, N\Big)$.
By the higher order interior and boundary regularity (see \cite{LW}), we conclude that
$u\in C^\infty\Big(\overline{B_\frac12^+}, N\Big)$. \hfill\qed  

\setcounter{section}{3} \setcounter{equation}{0}
\section{Proof of Theorem \ref{bdry_reg}}

In this section, we will apply Lemma 3.1 and proposition 2.1 to show that any stationary biharmonic
map satisfying (\ref{bdry_mono}) is smooth in $\overline{\Omega}$, away from a closed subset
with $(n-4)$-dimensional Hausdorff measure zero. 

Similar to the handling of interior regularity given by \cite{CWY} Lemma 4.8, \cite{W1} Lemma 5.3, and \cite{SM} appendix B, we need to establish
the following Morrey norm bound at the boundary.
\begin{lemma} Under the same assumptions as Theorem \ref{bdry_reg},
there exist $\epsilon_0>0$, $\theta_0\in (0,1)$, and $R_1=R_1(R_0,\epsilon_0)$ 
such that if for $x_0\in\partial\Omega$
and $0<R\le R_1$,
\begin{equation}\label{bdry_small3}
{R}^{4-n}\int_{\Omega\cap B_{R}(x_0)}\left(|\nabla^2 u|^2+r^{-2}|\nabla u|^2\right)
\le\epsilon_0^2,
\end{equation}
then 
\begin{equation}\label{bdry_small4}
\Big\|\nabla u\Big\|_{M^{4,4}(\Omega\cap B_{\theta_0 R}(x_0))}
+\Big\|\nabla^2 u\Big\|_{M^{2,4}(\Omega\cap B_{\theta_0 R}(x_0))}\le C\epsilon_0.
\end{equation}
\end{lemma}
\pf By the boundary monotonicity inequality (\ref{bdry_mono})
along with suitable translations and scalings, a crucial step to establish (\ref{bdry_small4}) is to obtain
certain control of $r^{4-n}\int_{\Omega\cap B_r(x_0)}|\nabla u|^2$. 

There are two different ways to show (\ref{bdry_small4}): the first one is similar to that by \cite{CWY} Lemma 4.8 and \cite{W1} Lemma 5.3, and the second one is similar to the new, simpler approach by \cite{SM} appendix B. Here we provide the second one, which is a slight modification of \cite{SM}. 

First, let's define 
\begin{equation} \label{sing_set}
\Sigma:=\Big\{
x\in\overline{\Omega}: \liminf_{r\downarrow 0} 
\ r^{4-n}\int_{\Omega\cap B_r(x_0)}(|\nabla^2 u|^2+r^{-2}|\nabla u|^2)>0
\Big\}.
\end{equation}
Then, it is a standard fact  that $H^{n-2}(\Sigma)=0$ (see Evans-Gariepy \cite{EG}).  
In particular, we have that for $H^{n-1}$ a.e. $x_0\in\partial\Omega$, 
\begin{equation}\label{small_bdry_pt}
\liminf_{r\downarrow 0} 
\ r^{4-n}\int_{\Omega\cap B_r(x_0)}(|\nabla^2 u|^2+r^{-2}|\nabla u|^2)=0. 
\end{equation}
With the help of the interior argument by \cite{SM} appendix B, it suffices to show 
that if $x_0\in\partial\Omega$ is such that (\ref{small_bdry_pt}) holds
and $R_0>0$ is such that (\ref{bdry_mono}) holds, then there exists $R_1>0$
depending on $R_0,\epsilon_0$, and $\phi$ such that  for any $0<r\le R$, 
there exist $\frac{r}2<r_0<r$ and $C=C(n,\Omega, N)>0$  such that
\begin{equation}\label{small_bdry_pt1}
 r_0^{4-n}\int_{\Omega\cap B_{r_0}(x_0)}|\nabla^2 u|^2
+r_0^{3-n}\int_{\partial B_{r_0}(x_0)\cap\Omega}|\nabla u|^2\le C\epsilon_0^2.
\end{equation}
In fact, since $u=\phi$ on $\partial\Omega\cap B_{r}(x_0)$, then by the $H^2$-estimate of
Laplace equation, one has
\begin{eqnarray}\label{bdry_h2_est}
&&r^{4-n}\int_{B_{\frac{r}3}(x_0)\cap\Omega}|\nabla^2 u|^2+r^{2-n}
\int_{B_{\frac{r}3}(x_0)\cap\Omega}|\nabla u|^2\nonumber\\
&\leq& C\Big[r_0^{4-n}\int_{B_{r_0}(x_0)\cap\Omega}|\nabla^2 u|^2+r_0^{2-n}
\int_{B_{r_0}(x_0)\cap\Omega}|\nabla u|^2+\|\phi\|_{C^2}^2 r_0^2\Big]
\le C\epsilon_0^2.
\end{eqnarray}
This, combined with Nirenberg's inequality, implies
\begin{equation}
r^{4-n}\int_{B_\frac{r}3(x_0)\cap\Omega}|\nabla u|^4
\le C\Big[r^{4-n}\int_{B_{\frac{r}3}(x_0)\cap\Omega}|\nabla^2 u|^2+r^{2-n}
\int_{B_\frac{r}3(x_0)\cap\Omega}|\nabla u|^2+\|\phi\|_{C^2}^2 r^2\Big]\le C\epsilon_0^2.
\end{equation}
Now we want to establish (\ref{small_bdry_pt}) by applying (\ref{bdry_mono}). For simplicity, assume
$x_0=0$ and set $\sigma(r)=\sigma_1(r)+\sigma_2(r)$, where
\begin{equation}
\sigma_1(r)=r^{4-n}\int_{B_r(0)\cap\Omega}|\nabla^2 (u-\phi)|^2
+r^{3-n}\int_{\partial B_r(0)\cap\Omega}|\nabla (u-\phi)|^2,
\end{equation}
and
\begin{equation}
\sigma_2(r)=r^{3-n}\int_{\partial B_r(0)\cap\Omega}
\left(2x^i(u-\phi)_{ij}(u-\phi)_j+3|\nabla (u-\phi)|^2-4r^{-2}|x^i(u-\phi)_i|^2\right).
\end{equation}
Then (\ref{bdry_mono}) implies that for any $0<r<R\le R_0$,
\begin{eqnarray}
\label{bdry_mono5}
&&\sigma(r)+\int_{\Omega\cap(B_R(0)\setminus B_r(0))} \left(\frac{|(u-\phi)_j+x^i (u-\phi)_{ij}|^2}{|x|^{n-2}}+(n-2)\frac{|x^i(u-\phi)_i|^2}{|x|^n}\right)
\nonumber\\
&&\leq CRe^{CR}+e^{CR}\sigma(R).
\end{eqnarray}
Since for a good radius $r>0$ we can bound
$$|\sigma(r)|
\le C[r^{4-n}\int_{\Omega\cap B_{2r}(0)}(|\nabla^2 u|^2+r^{-2}|\nabla u|^2)]+C\|\phi\|_{C^2}^2r^2,$$
(\ref{small_bdry_pt}) implies
$$\liminf_{r\downarrow 0}|\sigma(r)|=0.$$
Therefore (\ref{bdry_mono5}) and (\ref{bdry_small3}) imply that for $0<r\le R$,
\begin{eqnarray}
&&\int_{\Omega\cap B_{r}(0)}\left(\frac{|(u-\phi)_j+x^i (u-\phi)_{ij}|^2}{|x|^{n-2}}+(n-2)\frac{|x^i(u-\phi)_i|^2}{|x|^n}\right)\nonumber\\
&&\le C[\epsilon_0^2+\|\phi\|_{C^2}^2 r^2]\le C\epsilon_0^2.
\end{eqnarray}
Hence we obtain
\begin{equation}
\inf_{\frac{r}2\le\rho\le r}
\rho^{3-n}\int_{\partial B_\rho(0)\cap\Omega}
\left(\frac{|(u-\phi)_j+x^i (u-\phi)_{ij}|^2}{|x|^{n-2}}+(n-2)\frac{|x^i(u-\phi)_i|^2}{|x|^n}\right)
\le C\epsilon_0^2.
\end{equation}
Using (83) and (84) in \cite{SM} page 262, this implies
\begin{equation}
\sup_{\frac{r}2\le\rho\le r}\sigma_2(\rho)\ge -C\epsilon_0^2.
\end{equation}
This and the monotonicity inequality (\ref{bdry_mono5}) imply that there eixsts
$r_0\in [\frac{r}2, r]$ such that
\begin{equation}
\sigma_1(r_0)\le \sigma(R)+\sigma_2(r_0)\le C\epsilon_0^2.
\end{equation}
This implies (\ref{small_bdry_pt1}). \qed\\

\noindent{\bf Proof of Theorem \ref{bdry_reg}}: For $\epsilon_0>0$  given by Lemma 3.1, 
set
$$\mathcal S(u)
=\left\{x\in\overline\Omega: \liminf_{r\downarrow 0} r^{4-n}\int_{\Omega\cap B_r(x)}
(|\nabla^2 u|^2+|\nabla u|^4)\ge\epsilon_0^2\right\}.$$
Then it holds (cf. \cite{EG}) that $H^{n-4}(\mathcal S(u))=0$. It follows from the interior regularity
and the boundary regularity Lemma 3.1 that $u\in C^\infty(\overline\Omega\setminus
\mathcal S(u), N)$. \qed\\

\noindent{\bf Acknowledgement}. The third author is partially supported by NSF DMS 1000115. 
The work was completed when the first author visited University of Kentucky, partially supported by
NSF DMS 0601162. The first author would like to thank the department of Mathematics for its
hospitality.


\bigskip


\begin{thebibliography}{DU}

\bibitem{A} G. Angelsberg, {\em  A monotonicity formula for stationary biharmonic maps}.
Math. Z., 252 (2006), 287-293.

\bibitem{B} F. Bethuel, {\em On the singular set of stationary harmonic maps.}
Manus. Math. 78, no. 4 (1993): 417–43.

\bibitem{ADN} S. Agmon, A. Douglas, L. Nirenberg, 
{\em Estimates near the boundary for solutions of elliptic partial differential equations
satisfying general boundary conditions I}. Comm. Pure Appl. Math. 12 (1959) 623-727.

\bibitem{CWY} A. Chang, L. Wang, P. Yang, {\em A regularity theory of biharmonic maps}.
Comm. Pure Appl. Math. 52 (1999), no. 9, 1113-1137.

\bibitem{DS} A. Dall'Acqua, G. Sweers,
{\em  Estimates for Green function and Poisson kernels of higher-order 
Dirichlet boundary value problems.} 
JDE. Vol. 205, No. 2, 2004, 466-487.

\bibitem{E} L. C. Evans, {\em Partial regularity for stationary harmonic maps into spheres.}
Arch. Ration. Mech. Anal. 116, no. 2 (1991): 101–13.

\bibitem{EG} L. C. Evans, R. Gariepy, {Measure Theory and Fine Properties of Functions.}
Studies in Advanced Mathematics. Boca Raton, FL: CRC Press, 1992.

\bibitem{GS} A. Gastel, C. Scheven, {\em Regularity of polyharmonic maps in the critical dimension}. 
Comm. Anal. Geom. 17 (2) 2009, 185-226.

\bibitem{H} F. H\'elein, {\em R\'egularit\'e des applications faiblement harmoniques entre une surface et une vari\'et\'e riemannienne}.
C. R. Acad. Sci. Paris S\'er. I Math. 312 (1991), pp. 591–596.

\bibitem{Kra} J. Krasovskii, {\em Isolation of the singularity in Green's function}, Izv. Akad. Nauk SSSR Ser. Mat.
31 (1967) 977-1010 (in Russian).

\bibitem{K} Y. Ku, {\em Interior and boundary regularity of intrinsic biharmonic maps to spheres}.
Pacific J. Math. 234 (2008), pp. 43–67.

\bibitem{LR} T. Lamm, T. Rivier\'e, {\em Conservation laws for fourth order systems in four dimensions} Comm. PDE. 33 (2008), pp. 245–262.

\bibitem{LW} T. Lamm, C. Y. Wang, {\em Boundary regularity for polyharmonic maps in the critical dimension.}
Adv. Cal. Var. 2, 2009, No. 1, 1-16.

\bibitem{P} P. Price, {\em A monotonicity formula for Yang-Mills fields}. Manus. Math. 43 (1983),
131-166.

\bibitem{RT} T. Rivier\'e, {\em Conservation laws for conformally invariant variational problems.}
Invent. Math. 168 (1), (2007) 1-22.

\bibitem{RS} T. Rivier\'e, M. Struwe, {\em Partial regularity of stationary harmonic maps and related problems}. Comm. Pure Appl. Math. 61, no. 4 (2008), 451– 63.

\bibitem{SM} M. Struwe, {\em Partial regularity for biharmonic maps, revisited}.
C.V.\& P.D.E. 33 (2008), No. 2, 249-262.

\bibitem{Str} P. Strzelecki, {\em On biharmonic maps and their generalizations}. C.V.\& P.D.E.
18 (4), 401-432 (2003).

\bibitem{W} C. W. Wang, {\em Biharmonic maps from R4 into a Riemannian manifold.}
Math. Z. 247 (2004), pp. 65–87.

\bibitem{W1} C. Y. Wang,  {\em Stationray biharmonic Maps from
$\mathbb{R}^m$ into a Riemannian Manifold.}  Comm. Pure
Appl. Math. 57, (2004) 0419-0444.

\bibitem{W2} C. Y. Wang, {\em  Remarks on biharmonic maps into spheres.}
C.V.\& P. D.E., 21 (2004) 221-242.

\bibitem{W3} C. Y. Wang, {\em Boundary partial regularity for a class of harmonic maps.}
Comm. PDE 24 (1999), no. 1-2, 355-368.

\end{thebibliography}
 \end{document}